\documentclass{amsart}

\usepackage{amssymb}
\usepackage{graphics, latexsym, tabularx, shapepar, enumerate}
\usepackage[all]{xy} %\UseAllTwocells \SilentMatrices

\renewcommand{\P}{{\mathbb P}}
\newcommand{\Q}{{\mathbb Q}}

\newcommand{\Z}{{\mathbb Z}}

\newcommand{\Gal}{\mathrm{Gal}}

\newcommand{\Out}{\mathrm{Out}}

\newcommand{\Spec}{\mathrm{{\bf Spec}}\,}

\newtheorem{theorem}{Theorem}
\newtheorem{lemma}[theorem]{Lemma}
\newtheorem{corollary}[theorem]{Corollary}

\begin{document}

\title{On the torsion of Jacobians of principal modular curves of level $3^n$}
\author{Matthew Papanikolas}
\address{Department of Mathematics\\
Texas A{\&}M University\\
College Station, TX 77843}
\email{map@math.tamu.edu}
\author{Christopher Rasmussen}
\address{Department of Mathematics\\
Rice University\\
Houston, TX 77005}
\email{crasmus@math.rice.edu}

\begin{abstract}
We demonstrate that the $3$-power torsion points of the Jacobians of
the principal modular curves $X(3^n)$ are fixed by the kernel of the
canonical outer Galois representation of the pro-$3$ fundamental
group of the projective line minus three points. The proof proceeds
by demonstrating the curves in question satisfy a two-part criterion
given by Anderson and Ihara. Two proofs of the second part of the
criterion are provided; the first relies on a theorem of Shimura,
while the second uses the moduli interpretation.
\end{abstract}


\date{September 30, 2005}

\maketitle

\section{Introduction}

Let $\ell$ be a prime, and let $G_\Q$ denote the absolute Galois
group $\Gal (\bar{\Q}/\Q)$. We let $\pi_1$ denote the pro-$\ell$
algebraic fundamental group of the projective line minus three
points, $\pi_1 := \pi_1^\ell \bigl( \P^1_{\bar{\Q}} \smallsetminus
\{0, 1, \infty \} \bigr)$. This group may also be identified with
the Galois group $\Gal \bigl(M / \bar{\Q}(t) \bigr)$, where $M$ is
the maximal pro-$\ell$ extension of $\bar{\Q}(t)$ unramified away
from the places $t=0$, $t=1$, and $t = \infty$. The tower of Galois
extensions $\Q(t) \subseteq \bar{\Q}(t) \subseteq M$ yields an exact
sequence of profinite groups
\begin{equation*}
\xymatrix{ 1 \ar[r] & \pi_1 \ar[r] & \Gal \bigl( M / \Q(t) \bigr)
\ar[r] & G_\Q \ar[r] & 1. }
\end{equation*}
In the natural way, we define a Galois representation
\begin{equation*}
\rho_\ell \colon G_\Q \longrightarrow \Out (\pi_1).
\end{equation*}

The structure of $\rho_\ell$ is of active interest. Although it is not
injective, $\rho_\ell$ has a rich structure in the sense that the
fixed field of its kernel, denoted $\Omega_\ell$, is a quite large
subfield of $\bar{\Q}$. In particular, $\Omega_\ell$ is known to be an
infinite, non-abelian pro-$\ell$ extension of the field of
$\ell$-power roots of unity $\Q \bigl( \mu_{\ell^\infty} \! \bigr)$,
unramified away from $\ell$ \cite{Anderson:1988}. Let $\Lambda_\ell$
be the \emph{maximal} pro-$\ell$ extension of $\Q \bigl(
\mu_{\ell^\infty} \! \bigr)$ unramified away from $\ell$.  It is an
open problem, originally posed by Ihara \cite{Ihara:1986}, whether
$\Omega_\ell$ and $\Lambda_\ell$ coincide. Recent work of Sharifi
\cite{Sharifi:2002} has connected the answer to certain conjectures of
Deligne and Greenberg, but the problem itself remains open.

Another line of investigation for this question is the following:
which subfields of $\Lambda_\ell$ are also subfields of
$\Omega_\ell$? A standard recipe for large subfields of
$\Lambda_\ell$ is well-known. Let $C$ be a complete non-singular
curve defined over a number field $k \subseteq \Omega_\ell$, whose
Jacobian $J$ has good reduction away from the primes of $k$ above
$\ell$ (for example, take $C$ itself to have good reduction away
from $\ell$), and suppose further that
\[
  k \bigl( J[\ell] \bigr) \subseteq \Lambda_\ell.
\]
Then the maximality of $\Lambda_\ell$ and the work of Serre and Tate
\cite{Serre:1968} imply that the larger field $k \bigl(
J[\ell^\infty] \bigr)$ also lies in $\Lambda_\ell$, as it must be
pro-$\ell$ over $k \bigl( J[\ell] \bigr)$, and unramified away from
$\ell$.

Hence, any such curve $C$ provides a natural subfield to consider,
and in many cases the (possibly) stronger containment
\begin{equation}\label{eq:Omega}
k \bigl( J[\ell^\infty] \bigr) \subseteq \Omega_\ell
\end{equation}
is known. For example, Anderson and Ihara \cite{Anderson:1988}
verified the containment \eqref{eq:Omega} for Fermat curves and
Heisenberg curves of prime power level, as well as for modular
curves of level $2^n$. In \cite{Rasmussen:2004} the second author
has shown that \eqref{eq:Omega} holds for elliptic curves over $\Q$
when $\ell = 2$. The cases of elliptic curves over $\Omega_2$, and
of elliptic curves over $\Q$ when $\ell = 3$ are partially treated
in \cite{Rasmussen:Thesis}.

\subsubsection*{Geometric $\ell$-Covers}

Let $B$ denote $\P^1_\Q \smallsetminus \{0, 1, \infty \}$. Consider
the pairs $(C,g)$, where $C$ is a complete nonsingular curve,
defined over a field $k \subseteq \bar{\Q}$, and
\begin{equation*}
g \colon C \longrightarrow \P^1_k
\end{equation*}
is a morphism defined over $k$. We say such a pair is a geometric
$\ell$-cover if the base extension morphism
\begin{equation*}
\bar{g} := g \otimes_k \bar{\Q} \colon C_{\bar{\Q}} \longrightarrow
\P^1_{\bar{\Q}}
\end{equation*}
has Galois closure of $\ell$-power degree, and $\bar{g}$ ramifies
only over the set $\{0, 1, \infty \}$. This second condition is
equivalent to the requirement that $\bar{g}$ restricts to an \'etale
morphism over $B_{\bar{\Q}}$.

For all choices of $C$ for which \eqref{eq:Omega} is known, the
proof of \eqref{eq:Omega} proceeds by demonstrating the curve as a
geometric $\ell$-cover of the projective line minus three points.
Then, one uses the representation $\rho_\ell$, which relates the
arithmetic of such a curve (via Galois action) to its geometry as
such a cover. In fact, Anderson and Ihara provided a sufficient
criterion to determine when the existence of such an $\ell$-cover
implies \eqref{eq:Omega}. The following theorem, which is the main
result of this article, follows the same approach in its proof.

\begin{theorem}\label{thm:main}
Let $n \geq 1$. Let $X(3^n)$ be the principal modular curve of level
$3^n$, and let $J$ be its Jacobian variety. Then $J[3^\infty]$ is
rational over $\Omega_3$.
\end{theorem}

In this article, we let $X(N)$ denote the principal modular curve,
which is associated to the group $\Gamma(N)$ of $2 \times 2$ matrices
congruent to the identity modulo $N$ and is defined over the
field of $N$-th roots of unity $\Q(\zeta_N)$.  Similarly, we let
$X_0(N)$ denote the modular curve over $\Q$, which is associated to
the congruence subgroup $\Gamma_0(N)$ of matrices that are upper
triangular modulo $N$.

\section{Criterion of Anderson and Ihara}

Throughout we let $\Delta$ denote the field of Puiseux series over
$\Q$ in $1/t$:
\begin{equation*}
\Delta := \Q((1/t))[t^r : r \in \Q].
\end{equation*}
Let $\cdot$ indicate the operation of field compositum. The result
of Anderson and Ihara \cite{Anderson:1988} mentioned above is the
following.

\begin{theorem}[Anderson-Ihara \cite{Anderson:1988}, 3.8.1]\label{thm:AIC}
%%Let $\Delta$ denote the field of Puiseux series over $\Q$.
Let $k$ be a subfield of $\Omega_\ell$, and let $C$ be a smooth
projective curve over $k$, with Jacobian variety $J$. Suppose that
\begin{enumerate}[(i)]
\item There exists a morphism $g \colon C_{\Omega_\ell} \to
\P^1_{\Omega_\ell}$ unramified outside $\{0,1,\infty\}$, which is a
geometric $\ell$-cover.
\item There exists a point $y \in C(\Omega_\ell \cdot \Delta)$,
such that $g(y) = t \in \P^1(\bar{\Q} \cdot \Delta)$.
\end{enumerate}
Then the $\ell$-power torsion points of $J$ are rational over
$\Omega_\ell$.
\end{theorem}

\subsection{Remark on the $\ell = 2$ case}

It is worthwhile to review the approach of Anderson and Ihara in the
case of the principal modular curves of level $2^n$. For any fixed
$n$, there is a natural tower of degree $2$ covers, defined over $\Q
\bigl( \mu_{2^\infty}\! \bigr)$,
\begin{equation*}
\xymatrix{ X(2^n) \ar[r] & X(2^{n-1}) \ar[r] & \cdots \ar[r] & X(2)
},
\end{equation*}
which ramify only over the cusps of $X(2)$. But $X(2) \smallsetminus
\{\mbox{cusps} \} \cong \P^1 \smallsetminus \{0,1,\infty\}$ over $\Q$,
and so this gives a $2$-cover $X(2^n) \to \P^1$ over $\Omega_2$ (this
covering is itself Galois, since the congruence subgroup $\Gamma(2^n)$
is normal in $\mathrm{SL}_2(\Z)$). The second requirement of the
criterion now follows from the fact that the function field for $X(2^n)$
is generated by modular forms whose Fourier expansions have especially
nice properties. The details are very similar (and in fact simpler)
than the argument given in the following section for $3$-power level.

\subsection{The First Condition}

For the remainder of the article, we let $L = \Q \bigl(
\mu_{3^\infty}\! \bigr)$. Fix $n \geq 1$. At first glance, the above
approach for $\ell = 2$ seems unavailable for $\ell = 3$, as $X(3)$
possesses $4$ cusps. However, the modular curve $X_0(3)$ possesses
two cusps and an elliptic point, and so in fact the composition of
maps defined over $L$
\begin{equation*}
\xymatrix{ g \colon X(3^n) \ar[r] & X(3^{n-1}) \ar[r] & \cdots
\ar[r] & X(3) \ar[r] & X_0(3) \ar[r]^\simeq & \P^1 }
\end{equation*}
demonstrates the desired $3$-cover for $X(3^n)$. All the maps are
the natural ones; the final isomorphism will be a composition of the
Hauptmodul for $X_0(3)$ with a certain linear fractional
transformation.

The Hauptmodul for $X_0(3)$ that we will use is the function
\begin{equation} \label{eq:Haupt}
 h(z) := \biggl( \frac{\eta(z)}{\eta(3z)} \biggr)^{12},
\end{equation}
where $\eta(z)$ is the Dedekind eta function.  It is easily seen to
be a modular function for $\Gamma_0(3)$.  On $X_0(3)$ the function
$h$ has a simple zero at the cusp $0$ and a simple pole at the cusp
$\infty$ with residue $1$.  Moreover by checking Fourier expansions,
it is a straightforward matter to observe (as in \cite{Chen:1996})
that
\[
  j = \frac{(h+27)(h+243)^3}{h^3},
\]
where $j$ is the classical $j$-function.  Thus we have a commutative
diagram of curves over $\Q$,
\begin{equation}\label{comm_diagram}
\xymatrix{
& X_0(3) \ar[r]^h \ar[d]_{\mathbf{pr}} & \P^1_\Q \ar[d]^{x
\mapsto x^{-3}(x+27)(x+243)^3 } \\
& X(1) \ar[r]_j & \P^1_\Q }
\end{equation}
Our construction of the map $g$ above will be complete once we
compose $h$ with the appropriate linear fractional transformation
that takes $0$ to $0$, $\infty$ to $\infty$, and the elliptic point
of $X_0(3)$ to $1$.  The function $g$ will satisfy the first part of
Theorem \ref{thm:AIC} for $X_0(3^n)$, provided that this linear
fractional transformation is defined over $\Omega_3$.

Let $e$ denote the unique elliptic point on $X_0(3)$. We are
interested in where $e$ and the two cusps $0$ and $\infty$ are
mapped under $h$. The image of $e$ under $\mathbf{pr}$ must be
$\rho$, the elliptic point of order $3$ on $X(1)$.  Since $j(\rho) =
0$, $h(e)$ must be either $-27$ or $-243$. But $e$ is not a
ramification point of $\mathbf{pr}$, and so $h(e) = -27$.  Thus the
necessary fractional linear transformation is merely the scaling $x
\mapsto -\frac{1}{27}x$, which is defined over $\Q$. Hence, we
certainly have $g \colon X(3^n)_L \rightarrow \P^1_L$ defined over
$L \subseteq \Omega_3$. Thus, $g$ gives a geometric $3$-cover for
$X_0(3^n)$ satisfying the first condition of Theorem \ref{thm:AIC}.

\subsection{The Second Condition: Proof 1}

Here and in the next section, we prove that the covers just
constructed satisfy the second condition of Theorem \ref{thm:AIC}.
Certainly, it is enough to produce a morphism $y$ which makes the
following diagram commute:
\begin{equation*}
\xymatrix{
\Spec (\Omega_3 \cdot \Delta) \ar[d] \ar@{-->}[rr]^{\exists?\;y} &
& X(3^n)_{L(t)} \ar[d]^{g \otimes_L L(t)} \\
\Spec (L \cdot \Delta) \ar[rr]^t & & \P^1_{L(t)}
}
\end{equation*}
Let $z$ be an indeterminate, and identify $\Spec \bigl( L(t)[z^{-1}]
\bigr)$ as an affine neighborhood of $\infty \in \P^1_{L(t)}$. Then
the morphism $t$ in the above diagram corresponds to the ring
homomorphism given by evaluation:
\begin{equation*}
L(t)[z^{-1}] \longrightarrow \Omega_3 \cdot \Delta, \qquad \qquad z
\mapsto t.
\end{equation*}
Similarly, $\Spec \bigl( L(t)[h^{-1}] \bigr)$ can be identified with
an affine neighborhood of $\infty \in X_0(3)_{L(t)}$. The
restriction of the \emph{morphism} $h$ on this neighborhood
corresponds to
\begin{equation*}
L(t)[z^{-1}] \longrightarrow L(t)[h^{-1}], \qquad \qquad z \mapsto -
\tfrac{h}{27}.
\end{equation*}
Let $\theta$ be a generator for the function field of
$X(3^n)_{L(t)}$ as an extension of the function field for
$X_0(3)_{L(t)}$. Then we can give an affine open subscheme of
$X(3^n)_{L(t)}$ by $\Spec \bigl( L(t)[h^{-1}, \theta] \bigr)$, and
the natural map to $\Spec \bigl( L(t)[h^{-1}] \bigr)$ is a
restriction of the covering $X(3^n)_{L(t)} \rightarrow
X_0(3)_{L(t)}$ given in the previous section.

A theorem of Shimura \cite[6.9]{Shimura:1971} guarantees that the
series expansion for $\theta$ in fractional powers of $q$ lies in
the ring $\Omega_3 \bigl( \bigl( q^{1/3^n} \bigr) \bigr)$. In fact,
one has the stronger condition that the coefficients of this series
lie in the field $L$.

Hence, it is possible to give a ring homomorphism
\begin{equation*}
y \colon L(t)[h^{-1}, \theta] \longrightarrow \Omega_3 \cdot \Delta
\end{equation*}
as follows. First there is a homomorphism
\begin{equation*}
\beta \colon L(t)[h^{-1}, \theta] \longrightarrow \Omega_3 \bigl(
\bigl( q^{1/3^n} \bigr) \bigr)
\end{equation*}
given by mapping $h$ and $\theta$ to their Fourier series in $q$,
and sending $t$ to the Fourier series for $-\frac{h}{27}$.

By \eqref{eq:Haupt} the Fourier series for the Hauptmodul $h$ has
the form
\begin{equation*}
h = h(q) = \frac{1}{q} - 12 + 54q - 76q^2 - 243q^3 + 1188q^4 +
\cdots \in q^{-1}\Q[[q]],
\end{equation*}
and so there is an equivalent relation giving $q$ as a $\Q$-rational
series
\begin{equation*}
q = h^{-1} - 12h^{-2} + 198h^{-3} - 3748h^{-4} + 76629h^{-5} +
\cdots \in \Q[[h^{-1}]].
\end{equation*}
Hence, we can write
\begin{equation*}
q^{1/3^m} = h^{-1/3^m} \sum_{i \geq 0} b_{i,m} h^{-i} \in h^{-1/3^m}
\cdot \Q[[h^{-1}]].
\end{equation*}
Define a homomorphism
\begin{equation*}
\alpha \colon \Omega_3 \bigl( \bigl( q^{1/3^n} \bigr) \bigr)
\longrightarrow \Omega_3 \cdot \Delta
\end{equation*}
by sending $q^{1/3^n}$ to its $h$-series evaluated at $h^{1/3^n} =
 -27^{1/3^n} t^{1/3^n}$. The composition $\alpha \circ \beta$ gives the
homomorphism $y$. By abuse of notation, we use ring homomorphisms to
denote the associated morphisms of schemes. We have a commutative
diagram:
\begin{equation*}
\xymatrix @C=0.3in{
& & \Spec \left( \Omega_3 \bigl( \bigl( q^{1/3^n} \bigr) \bigr)
\right) \ar[d]^\beta & & \\
\Spec (\Omega_3 \cdot \Delta) \ar[rru]^\alpha \ar[d] \ar@{-->}[rr]^y
& & \Spec \left( L(t)[h^{-1}, \theta ] \right) \ar[d] & \subseteq
& X(3^n)_{L(t)} \ar[d] \ar@/^2.5pc/[dd]^{g \otimes_L L(t)} \\
\Spec (L \cdot \Delta) \ar@{=}[d] \ar[rr]_{h \mapsto -27t} & & \Spec
\left( L(t) [ h^{-1} ] \right) \ar[d]^{z \mapsto - \frac{h}{27}} &
\subseteq & X_0(3)_{L(t)} \ar[d] \\
\Spec (L \cdot \Delta) \ar[rr]_{z \mapsto t}& & \Spec \left(
L(t)[z^{-1}] \right) & \subseteq & \P^1_{L(t)} \\
}
\end{equation*}
Of course, this $y$ now has the required properties, and so the
second part of the criterion is satisfied. Hence, Theorem
\ref{thm:main} follows immediately.

\subsection{The Second Criterion: Proof 2}

The previous proof follows closely the original approach of Anderson
and Ihara in demonstrating the analogous result for the curves
$X(2^n)$. In what follows, we present an alternative proof.
Essentially, the dependence on the result of Shimura is removed by
using the moduli interpretation for the curves in question.

Fix $n \geq 1$. We would like to find points in $X(3^n)$ that map to
$t$ under our cover $g$. Points in $X(3^n)$ are represented by pairs
$(E, \{P, Q\})$, where $E$ is an elliptic curve, and $P, Q \in
E[3^n]$ provide a basis for the $3^n$-torsion of $E$. Similarly, the
points in $X_0(3)$ are represented by pairs $(E,C)$, where $C$ is a
subgroup of order $3$ in $E$. The natural maps in our cover $g$ are:
\begin{equation*}
\begin{split}
X \left( 3^n \right) \rightarrow X \left( 3^{n-1} \right) \qquad &
\qquad \bigl[ E, \{P, Q \} \bigr] \mapsto \bigl[ E, \{ [3]P, [3]Q \} \bigr], \\
X(3) \rightarrow X_0(3) \qquad & \qquad \bigl[ E, \{P', Q' \} \bigr]
\mapsto \bigl[ E, \langle P' \rangle \bigr].
\end{split}
\end{equation*}
As each of these maps is applied to a point $x$ on the modular
curve, the elliptic curve representing $x$ does not (have to)
change. Hence, we can find one elliptic curve $E$ which will,
together with various choices of torsion points on $E$, represent
all points on the fiber of $g$ over $t$. Such an $E$, together with
some cyclic subgroup $C$ of order $3$, must satisfy $h([E,C]) =
-27t$. Hence, $j(E) = f(-27t)$, where $f$ is the rational map on the
right side of the diagram \eqref{comm_diagram}. Hence,
\begin{equation*}
j(E) = \frac{(-27t +27)(-27t + 243)^3}{(-27t)^3} = -27
\frac{(t-1)(t-9)^3}{t^3}.
\end{equation*}
It is easy enough to write down an elliptic curve $E$ with such a
$j$-invariant using the Deuring normal form
\cite[Ex.~3.23]{Silverman:1986}. We choose
\begin{equation*}
E \colon y^2 + \alpha x y + y = x^3, \qquad \alpha = 3 \left(1 -
\tfrac{1}{t} \right)^{1/3}.
\end{equation*}
Let $\pi = t^{-1}$. Note that $\alpha \in R = \Omega_3[[\pi]]$,
since as a formal series
\begin{equation*}
\alpha = 3 - \pi - \frac{1}{3}\pi^2 - \frac{10}{9}\pi^3 - \cdots.
\end{equation*}
Hence, we may think of $E$ as an elliptic curve over the Laurent
series field $\Omega_3((\pi))$, which is a complete discrete
valuation ring with respect to the $\pi$-adic valuation. It is
simple to verify that the curve $E$ has split multiplicative
reduction at $\pi$.

Hence, the points on $X(3^n)$ which are mapped to $t \in \P^1(\bar{\Q}
\cdot \Delta)$ by $g$ are among the points given by the classes
$\mathcal{E} = \bigl[ E, \{P, Q \} \bigr]$, where $\{P, Q\}$ is any
basis for $E[3^n]$. We would like to show that at least one such
$\mathcal{E}$ lies in $X(3^n)(\Omega_3 \cdot \Delta)$. To do so, we
first observe that $E$ is already defined over $\Omega_3((\pi))
\subseteq \Omega_3\cdot\Delta$.  Thus it is enough to show that there
exists a basis $\{P, Q\}$ with $P, Q \in E( \Omega_3 \cdot
\Delta)$. Clearly, $P$ and $Q$ lie in $E (\bar{\Q} \cdot
\Delta)$. Hence, we can show they lie in $E(\Omega_3 \cdot \Delta)$ by
checking invariance under the Galois action of elements $\sigma \in
\Gal(\bar{\Q}/\Omega_3)$.  We will actually prove the following
stronger result, which will complete the proof of Theorem
\ref{thm:main}.
\begin{lemma}
Any point $P \in E[3^\infty]$ is defined over $\Omega_3 \cdot
\Delta$.
\end{lemma}
\begin{proof}
Let $K = \Omega_3((\pi))$ and $R = \Omega_3[[\pi]]$.  By Tate's
theorem on uniformization over complete discrete valuation fields (see
\cite[C.14.1]{Silverman:1986}), we know there is an isomorphism
\begin{equation*}
\varphi \colon \bar{K}^\times / q^\Z \rightarrow E(\bar{K})
\end{equation*}
of Galois modules, where $q \in \bar{K}^\times$ is an element such
that $j(q) = j(E)$. Let $A := \bar{K}^\times / q^\Z$. It is easy to
describe the preimage of the $m$-torsion of $E$ under $\varphi$; it
is given by points of the form $\eta^a q^{b/m}$, where $\eta$ is a
primitive $m$th root of unity, $q^{1/m}$ is a fixed $m$th root of
$q$ in $\bar{K}$, and $a$, $b \in (\Z / m\Z)^\times$.

It is a simple matter to solve for $q \in \bar{K}^\times$ from the
usual expansion
\begin{equation*}
j(q) = \frac{1}{q} + 744 + 196884q + \cdots.
\end{equation*}
In fact,
\begin{equation*}
q = -\frac{1}{27}\pi - \frac{4}{243}\pi^2 + \cdots = -\frac{1}{27}
\pi u, \qquad u \in 1 + \pi R.
\end{equation*}
Further, it is clear that the coefficients in the $\pi$-series for $q$
must be in $\Q$, because the same is true for both the
$q$-series of $j$, and the polynomial $j(E)$ in powers of $\pi$. But
in the form $q = - \frac{1}{27} \pi u$, we see that $E[3^n]$ is
isomorphic as a Galois module to the set
\begin{equation*}
A[3^n] = \left\{ \eta^a q^{b/3^n} \right. \left| \eta^{3^n} = 1;\,
a,b \in \Z/3^n\Z \right\}.
\end{equation*}
Now $\eta \in \Omega_3$, and $u \in 1 + \pi R$. Hence, roots of $u$
are also in $1 + \pi R$. Since $\pi^{1/3^n} = t^{-1/3^n}$ belongs to
the Puiseux series field over $\Omega_3$, it follows that all
$3^n$-th roots of $q$ are $(\Omega_3 \cdot \Delta)$-rational. Hence,
the $3^n$-torsion points in $E$ must be rational over this field
also.
\end{proof}

Again, this demonstrates the second part of the Anderson-Ihara
criterion. However, it should be noted that this result is stronger
than that of the preceding section. Previously, we demonstrated the
existence of a $(\Omega_3 \cdot \Delta)$-valued point $y$ in the
fiber over $t$. Here, we have demonstrated a large collection of
such points, at least one of which must lie in the fiber.

\section{Remarks}

\subsection{Proof 2 when $\ell = 2$}

Both of the methods presented for verifying the second part of the
Anderson-Ihara criterion also work for the case $\ell = 2$. Indeed,
the first version of the proof is quite direct in this case, because
there is no Hauptmodul involved in the covering morphism. This is
the approach used by Anderson and Ihara to originally treat the case
$X(2^n)$.

The second version of the proof may also be used for the $\ell = 2$
case. However, one must extend the covers down to $X_0(2)$ to use
the known relation between the Hauptmodul and the $j$-invariant. We
omit the details, but the appropriate elliptic curve has a model
given by
\begin{equation*}
E_t \colon \qquad y^2 = x(x-1)(x-\lambda), \qquad \lambda =
\frac{1}{2} \left( 1 + \sqrt{1 - t^4} \right).
\end{equation*}
Tate's theorem applies exactly as before to conclude the proof.

\subsection{Compatibility with covers}

Let $\varphi \colon X \to C$ be a morphism of complete non-singular
curves over $\Omega_\ell$, and let $J_X$ and $J_C$ denote their
respective Jacobians.  Suppose that the $\ell$-adic Tate module
$T_\ell(J_X)$ is known to be rational over the field $\Omega_\ell$.
Because the morphism $\varphi$ induces a natural inclusion of Galois
modules $T_\ell(J_C) \rightarrow T_\ell(J_X)$, it follows also that
$T_\ell(J_C)$ is rational over $\Omega_\ell$. This is equivalent to
saying that $J_C$ has $\Omega_\ell$-rational $\ell$-power torsion.

We can immediately apply this to the natural coverings of modular
curves
\begin{equation*}
X(\ell^n) \longrightarrow X_1(\ell^n) \longrightarrow X_0(\ell^n).
\end{equation*}
Hence, for $\ell = 2$ or $3$, the curves $X_0(\ell^n)$,
$X_1(\ell^n)$ all have Jacobians with $\Omega_\ell$-rational
$\ell$-power torsion.

Additionally, since every elliptic curve over $\Q$ is known to be
modular, there is always a morphism $X_0(\ell^n) \to E$ for any
elliptic curve $E$ over $\Q$ with good reduction away from $\ell$.
As a result, we have
\begin{corollary}
Let $E$ be an elliptic curve defined over $\Q$, with good reduction
away from $3$. Then $\Q \bigl( E[3^\infty] \bigr) \subseteq
\Omega_3$.
\end{corollary}
This improves the result of \cite{Rasmussen:Thesis}, which only
demonstrated this result for curves of $j$-invariant $0$.
Additionally, the authors have verified that these additional elliptic
curves themselves appear in the pro-$3$ tower of \'etale covers of the
projective line minus three points.

As far as the authors are aware, for every complete non-singular
curve $C$ defined over $\Omega_\ell$ for which the containment
\eqref{eq:Omega} holds, $C$ is known to be an $\ell$-cover of the
projective line minus three points. Several natural questions remain
unanswered and are worthy of further study: Are these two conditions
equivalent? If not, does the Anderson-Ihara Criterion exactly
describe the difference? Are there any interesting counter-examples?

\section*{Acknowledgements}

The authors would like to thank Greg Anderson, Ahmad El-Guindy, and
Yasutaka Ihara for their helpful conversations on the subject
matter.

%\bibliography{PRBib}
%\bibliographystyle{amsplain}

\providecommand{\bysame}{\leavevmode\hbox
to3em{\hrulefill}\thinspace}
\providecommand{\MR}{\relax\ifhmode\unskip\space\fi MR }
\providecommand{\MRhref}[2]{%
  \href{http://www.ams.org/mathscinet-getitem?mr=#1}{#2}
} \providecommand{\href}[2]{#2}

\end{document}